\documentclass{article}
\usepackage{amsmath,amsthm,amsfonts,graphicx}

\theoremstyle{definition}

\usepackage{multirow}
\begin{document}
\title{Fast CUR Approximation of Average Matrix and Extensions\thanks{This report is an outline of our paper \cite{PSZa}.}}
\author{Victor Y. Pan} 

\author{Victor Y. Pan$^{[1, 2],[a]}$,
John Svadlenka$^{[2],[b]}$, and  Liang Zhao$^{[2],[c]}$
\\
\and\\
$^{[1]}$ Department of Mathematics and Computer Science \\
Lehman College of the City University of New York \\
Bronx, NY 10468 USA \\
$^{[2]}$ Ph.D. Programs in Mathematics  and Computer Science \\
The Graduate Center of the City University of New York \\
New York, NY 10036 USA \\
$^{[a]}$ victor.pan@lehman.cuny.edu \\
http://comet.lehman.cuny.edu/vpan/  \\
$^{[b]}$ jsvadlenka@gradcenter.cuny.edu \\
$^{[c]}$  lzhao1@gc.cuny.edu 
} 
\date{}
\maketitle


\begin{abstract} 
CUR and low-rank approximations 
are among most fundamental subjects of  numerical 
linear algebra, with a
wide range of applications to a variety 
of highly important areas of modern computing,
which range from the machine
learning theory and neural networks to
 data mining and analysis. We first dramatically accelerate
 computation of such approximations for
the average input matrix, then show some narrow  classes of hard inputs
for our algorithms, and finally point out a tentative direction 
to narrowing such classes further
by means of  preprocessing with quasi Gaussian structured  multipliers.
Our extensive numerical tests with a variety of real world inputs
for regularization from Singular Matrix Database 
have consistently
produced reasonably close CUR approximations at a low computational cost.
There is a variety of efficient applications of our results and our techniques 
to important subjects of matrix computations. Our study provides new insights and 
enable design of faster algorithms for low-rank approximation by means of
 sampling and oversampling, for Gaussian elimination with no pivoting and block
 Gaussian elimination, and for the approximation of
railing singular spaces associated with the $\nu$ smallest singular values 
of a matrix having numerical nullity $\nu$.
 We conclude the paper with  novel extensions of our acceleration
to the Fast Multipole Method and the Conjugate Gradient Algorithms.
(This report is an outline of our paper \cite{PSZa}.)

 
\end{abstract}





 
\section{Background: CUR and low-rank approximation}\label{sbcgr}


{\em Low-rank approximation} of an $m\times n$ matrix $W$
having a small numerical rank $r$, that is, having
 a well-conditioned rank-$r$ matrix nearby, 
 is one of the most fundamental problems of 
numerical linear algebra \cite{HMT11}
 with a variety of  
  applications to highly important
areas of modern computing, which range from the machine
learning theory and neural networks \cite{DZBLCF14}, \cite{JVZ14}
to  numerous problems of data mining and analysis \cite{M11}. 

One of the most studied approaches to the solution of this problem
is given by  $CUR$ {\em approximation} 
where $C$ and $R$ are a pair of $m\times l$ 
 and $k\times n$ submatrices formed by $l$ columns and $k$ rows 
of the matrix $W$, respectively, and  $U$ is a $k\times l$
matrix such that $W\approx CUR$.

Every low-rank approximation allows very fast 
approximate multiplication 
of the matrix $W$ by a vector, but
CUR approximation is 
 particularly transparent and memory efficient.

The algorithms for computing it are characterized by the two main parameters:

(i) their complexity and

(ii) bounds on the error norms of the approximation.
  
We assume that  
$r\ll \min\{m,n\}$, that is, the integer $r$ is much smaller than $\min\{m,n\}$,
and we seek algorithms that use $o(mn)$ flops, that is,
much fewer than the information lower bound $mn$.


\section{State of the art and our progress}\label{ssartpr}


The algorithms of \cite{GE96}
and  \cite{P00} compute CUR approximations by using
 order of $mn\min\{m,n\}$ flops.\footnote{Here and hereafter
{\em ``flop"} stands for ``floating point arithmetic operation".}
 \cite{BW14} do this 
in $O(mn\log(mn))$ flops by using randomization. 

These are record upper bounds for computing a  CUR approximation
to {\em any input matrix} $W$, but
the user may be quite happy with having 
a close CUR approximations to {\em many matrices} $W$
that make up the class of his/her interest.
The information lower bound $mn/2$ 
(a flop involves at most two entries) does not apply  
to such a restricted input classes,
and we go well below it in our paper  \cite{PSZa}
(we must refer to that paper for technical details because 
of the limitation on the size of this submission).

We first formalize the problem of   CUR approximation
of an average $m\times n$ matrix of numerical rank $r\ll \min\{m,n\}$,
assuming the customary Gaussian (normal) probability distribution for 
its $(m+n)r$ i.i.d. input parameters.

Next we consider a two-stage approach:
(i) first fix a pair of integers $k\le m$ and $l\le n$
and compute a CUR approximation (by using the algorithms of 
\cite{GE96} or \cite{P00})
to a random $k\times l$ submatrix and then 
(ii) extend it to computing a
CUR approximation of an input matrix $W$ itself.
 
We must keep the complexity of Stage (i) low
and must extend
 the  CUR approximation from the submatrix 
 to the matrix $W$.
We prove that for a specific class of input matrices $W$ 
these two tasks are in conflict
(see Example 11 of \cite{PSZa}), 
but such a class of hard inputs is narrow,
because we prove that our algorithm 
produces a close approximation
to the average  $m\times n$ input matrix $W$ 
having numerical rank $r$. (We define such an
average matrix by assuming the standard Gaussian (normal)
probability distribution.)

By extending our two-stage
algorithms with the technique of \cite{GOSTZ10},
which we call {\em cross-approximation},
 we a little narrow the class of hard inputs
of Example 11 of \cite{PSZa} to the smaller class of  Example 14 of \cite{PSZa}
and moreover deduce a sharper bounds on the error of approximation by 
maximizing the {\em volume}
of an auxiliary $k\times l$ submatrix that defines a CUR approximation

In our extensive tests with a variety of real world input data  
for regularization of matrices from Singular Matrix Database,
 our fast algorithms consistently 
produce close CUR approximation.

Since our fast algorithms produce reasonably accurate
CUR approximation to the average input matrix,
the class of hard input matrices for these algorithms 
must be narrow, and we  
studied a tentative direction
towards further narrowing this input
class.
We prove that the algorithms are expected 
to output a close CUR approximation to any matrix $W$ 
if we pre-process it by applying Gaussian multipliers.
This is a nontrivial result 
of independent interest (proven on more than three pages), but
its formal support is only for application of Gaussian multipliers,
which is quite costly.

We hope, however, that we can still substantially narrow 
the class of hard inputs 
even if we replace
 Gaussian multipliers with
the products of  reasonable numbers 
of random bidiagonal matrices
and if we partly curb the permutation of these matrices.
If we achieve this, then preprocessing 
would become non-costly. 
This direction seems to be quite promising,
but still requires further work.

Finally, our algorithms can be extended to the
acceleration of various computational problems that are  known 
to have links to low-rank approximation, 
but in our concluding Section \ref{scncl} we describe a novel and rather unexpected extensions
to the acceleration of the Fast Multipole Method
and Conjugate Gradient Algorithms,\footnote{Hereafter
we use the acronyms FMM and CG.} both being among the most celebrated achievements
of the 20th century in Numerical Linear Algebra.


\subsection{Some related results on  matrix algorithms and our progress on 
other fundamental subjects of matrix computations}\label{srltwr}

 
A huge bibliography on CUR and low-rank approximation,
including the known best algorithms, which we already cited,
can be accessed from the papers \cite{HMT11}, \cite{M11}, \cite{BW14} and  \cite{W14}. 
 
Our main contribution is  dramatic acceleration of the  known algorithms.

Some of our techniques extend the ones of  \cite{PZ16}, \cite{PZ17}, and \cite{PZa},
where we also show
 duality of randomization and derandomization
and apply it to fundamental matrix computations.

In  \cite{PZ16} we prove  
that  preprocessing
with almost any well-conditioned multiplier of full rank
is 
as efficient on the average 
 for low-rank approximation as 
 preprocessing with a Gaussian one,
and then we propose  some new 
highly efficient sparse and structured multipliers.
Besides providing a new insight into the subject,
this motivates  the design of more efficient algorithms
and shows specific direction to this goal.

 We obtain similar progress in \cite{PZa} for
and  \cite{PZ17}
for  
preprocessing  Gaussian elimination with no pivoting and 
block Gaussian elimination.
We recall that Gaussian elimination with partial pivoting
is performed millions time per day, where 
pivoting, required for numerical stabilization, is 
frequently a bottleneck because
it interrupts the stream of arithmetic operations 
with foreign operations of comparison,
involves book-keeping, compromises data locality, and
increases communication overhead and data dependence. 

Randomized preprocessing is a natural 
substitution for pivoting, and  in \cite{PZa}
we show that 
Gaussian elimination with no pivoting as well as  
block Gaussian elimination (which is another valuable
algorithm and which also requires protection against numerical problems)
 are efficient on the average input
with preprocessing by any nonsingular and well-conditioned 
multipliers.  

 \cite{PZ17} obtains similar progress for
the  important subject of
the approximation of  trailing singular spaces 
associated with the $\nu$ smallest singular values 
of a matrix having numerical nullity $\nu$.

Our current progress greatly supersedes these earlier results, however, 
in terms of the scale of the acceleration of the known algorithms.

Our technique of representing random Gaussian multipliers 
as a product of random bidiagonal factors, 
our extension of CUR approximation to FMM and CG algorithms,
and our analysis of CUR approximation for the average input are new
and can have some independent interest.


\section{Conclusions}\label{scncl}


We dramatically accelerated  
the known algorithms for the fundamental problems of CUR 
and low-rank approximation 
in the case of the average input matrix and then
pointed out a  direction towards
heuristic extension of the resulting fast 
 algorithm  to a wider class of inputs
by applying quasi Gaussian preprocessing.
Our extensive tests for benchmark matrices 
of discretized PDEs have consistently supported
the results of our formal analysis.

Our study can be extended to a variety of important subjects of matrix computations.
Some of such extensions have been developed in papers   \cite{PZ16},  \cite{PZ17}
and  \cite{PZa}, and there are various challenging directions for further progress. 
In particular our accelerated CUR and low-rank approximation 
enables faster solution of some new important computational problems, thus
extending the long list of the known applications.
In the concluding section of  \cite{PSZa}, we  add two new highly
 important subjects to this long list.



\begin{thebibliography}{\hspace{1cm}}


\bibitem[BW14]{BW14}
C. Boutsidis, D. Woodruff,
Optimal CUR Matrix Decompositions,
{\em ACM Symposium on Theory of Computing (STOC 2014)},
353--362, ACM Press, New York, 2014.
Full version in arXiv:1405.7910v2  [cs.DS]  16 July 2014.


\bibitem[CI94]{CI94}
S. Chandrasekaran, I. Ipsen, 
On rank revealing QR factorizations, 
{\em SIAM J. Matrix Anal. Appl.}, {\bf 15}, 592--622, 1994.


\bibitem[DZBLCF14]{DZBLCF14}
E. Denton, W. Zaremba, J. Bruna, 
Y. LeCun, R. Fergus,
Exploiting Linear Structure Within Convolutional
Networks for Efficient Evaluation, {\em Proceedings of the 27th International
 Conference on Neural Information Processing Systems (NIPS'14)},
1269--1277, MIT Press, Cambridge, MA, USA, 2014.


\bibitem[GE96]{GE96}
M. Gu, S.C. Eisenstat, 
An efficient algorithm for computing a strong rank revealing QR factorization, 
{\em SIAM J. Sci. Comput.}, {\bf 17}, 848--869, 1996.


\bibitem[GOSTZ10]{GOSTZ10}
S. Goreinov, I. Oseledets, D. Savostyanov, E. Tyrtyshnikov, and N. Zamarashkin,
How to find a good submatrix, in 
{\em Matrix Methods: Theory, Algorithms, Applications}, 
(dedicated to the Memory of Gene Golub, edited by V. Olshevsky and E. Tyrtyshnikov), 
pages 247--256,
World Scientific Publishing, New Jersey, ISBN-13 978-981-283-601-4, ISBN-10-981-283-601-2 
 2010.


\bibitem[HMT11]{HMT11}
N. Halko, P. G. Martinsson, J. A. Tropp,
Finding Structure with Randomness: Probabilistic Algorithms
for Constructing
 Approximate Matrix Decompositions, 
{\em SIAM Review}, {\bf 53,~2}, 217--288, 2011.


\bibitem[JVZ14]{JVZ14}
M. Jaderberg, A. Vedaldi, A. Zisserman,
Speeding up Convolutional Neural Networks with Low Rank Expansions,
arXiv:1405.3866 CS, May 2014, and {\em BMVC} 2014.


\bibitem[M11]{M11}
M. W. Mahoney,
Randomized Algorithms for Matrices and Data,  
{\em Foundations and Trends in Machine Learning}, NOW Publishers, {\bf 3,~2}, 2011.
(Abridged version in: Advances in Machine Learning and Data Mining for Astronomy, 
edited by M. J. Way, et al., pp. 647-672, 2012.) 


\bibitem[P00]{P00}
C.-T. Pan,
On the existence and computation of
rank-revealing LU factorizations,
{\em Linear Algebra and its Applications}, {\bf 316}, 199--222, 2000.


\bibitem[PSZa]{PSZa}
V. Y. Pan, J. Svadlenka, L. Zhao,
Fast CUR Approximation of the Average Matrix and Extensions,
arXiv 1611.01391 (submitted November 3, 2016). 


\bibitem[PZ16]{PZ16}
V. Y. Pan, L. Zhao,
Low-rank Approximation of a Matrix:  
Novel Insights, New Progress, and Extensions,
Proc. of the {\em Eleventh International Computer Science Symposium in Russia  (CSR'2016)}, (Alexander Kulikov and Gerhard Woeginger, editors), St. Petersburg, Russia, June 2016,  {\em Lecture Notes in Computer Science (LNCS)}, {\bf 9691}, 352--366, Springer International Publishing, Switzerland (2016).  Also  arXiv:1510.06142 [math.NA], 
submitted on 21 Oct 2015, revised June 2016.


\bibitem[PZ17]{PZ17}
V. Y. Pan, L. Zhao,
New Studies of Randomized Augmentation and Additive Preprocessing,
{\em  Linear Algebra Appl.}, 2017, http://dx.doi.org/10.1016/j.laa.2016.09.035


\bibitem[PZa]{PZa}
V. Y. Pan, L. Zhao,
Numerically Safe Gaussian Elimination
with No Pivoting,
arxiv 1501.05385 CS, 
submitted on January 22, 2015, 
revised in June 2016.


\bibitem[W14]{W14}
D.P. Woodruff,
 Sketching As a Tool for Numerical Linear Algebra, 
{\em Foundations and Trends in Theoretical Computer Science}, 
{\bf 10,~1--2}, 1--157, 2014.


\end{thebibliography}
\end{document}